\documentclass[11pt]{amsart}
\usepackage{amsfonts,amssymb,amsmath}
\usepackage[all]{xy}
\begin{document}

\newcommand{\F}{{\mathbb F}}
\newcommand{\bF}{{\overline {\mathbb F}}}
\newcommand{\Hom}{{\rm Hom}}
\newcommand{\Q}{{\mathbb Q}}
\newcommand{\tL}{{\widetilde{\sf L}}}
\newcommand{\Gal}{{\rm Gal}}
\newcommand{\Symm}{{\rm Symm}}
\newcommand{\U}{{\rm U}}
\newcommand{\sL}{{\sf L}}
\newcommand{\x}{{\boldsymbol{\zeta}}} 
\newcommand{\Z}{{\mathbb Z}}

\title{A remark on Hopkins' chromatic splitting conjecture}
\author{Jack Morava}
\address{The Johns Hopkins University,
Baltimore, Maryland 21218}
\email{jack@math.jhu.edu}

\date{27 March 2015}
\begin{abstract}{Ravenel [4 \S 8.10] proved the remarkable fact that
the $K$-theoretic localization $L_K S^0$ of the sphere spectrum has 
\[
\pi_{-2} L_K S^0 \cong \Q/\Z \;.
\]
Hopkins' chromatic splitting conjecture [2] implies, more generally, that there 
are $3^{n-1}$ copies of $(\Q/\Z)_p$ in the homotopy groups of the $E(n)$-localization
of $S^0$; but where these copies occur can be confusing. We try here to simplify
this book-keeping.}\end{abstract} \bigskip

\thanks{\noindent This work was motivated by Hans-Werner Henn's talk on [1] at the August 2011 
Hamburg conference on structured ringspectra. I would also like to thank Agn\`es Beaudry
for helpful conversation.}

\maketitle \bigskip

\section{Introduction}

\noindent
This document is a footnote to Mark Hovey's account of Mike Hopkins' chromatic
splitting conjecture, but some of the notation here differs slightly from his:
\bigskip

\noindent
I'll write $[X,Y]$ for the spectrum of maps from $X$ to $Y$, and $[X,Y]_s$ 
for its $s$th homotopy group; $SX$ is the suspension of $X$ . The Bousfield 
localization of the sphere spectrum $S^0$ with respect to $E(n)$ will be our main
concern (the relevant prime $p$ will be suppressed) and for a general spectrum $E$ 
I'll simplify $L_E S^0$ to $L_E$ and write $L_n$ for $L_{E(n)}$. Since $E(n)$ is smashing, 
$L_n X = L_n \wedge X$ is a homology theory.\bigskip

\noindent 
There are standard homomorphisms $L_n \to L_{n-k}$ of ringspectra for $n \geq k$, 
and a fiber product square

\[
\xymatrix{
L_n \ar[d] \ar[r] & L_{K(n)} \ar[d] \\
L_{n-1} \ar[r] & L_{n-1}L_{K(n)} }
\]

\noindent
as well as a cofiber sequence [2 \S 4.1]

\[
\xymatrix{
S^{-1}L_{K(n)} \ar[r] & [L_{n-1},L_n] \ar[r] & L_n \;.}
\] 
\bigskip

\noindent
The splitting conjecture [H \S 4.2] asserts \medskip

$\bullet \;$ the existence of a monomorphism of an exterior algebra 

\[
E^*_{W(\F_p)}(\zeta_{2i+1} \:|\: 0 \leq i \leq n-1) \to H_c^*(S_n,W(\F_p))
\]

\noindent
[cf [3 \S 2.2.5]; $\zeta_i$ corresponds to Hovey's $x_{i+1}$] which becomes an isomorphism 
after rationalization, such that the products

\[
\zeta_I \; := \: \zeta_{2i_1+1} \cdots \zeta_{2i_l+1} 
\]

\noindent
(indexed by sequences $I = i_1,\dots,i_l$ with $0 \leq i_1 < \dots < i_l < n)$
survive the descent spectral sequence 

\[
H^*_c(S_n,W(\bF_p))^{\Gal(\bF_p/\F_p)} \Rightarrow \pi_*L_{K(n)}
\]

\noindent
to define maps

\[
\zeta_I : S^{-|I|} \to L_{K(n)} 
\]

\noindent
($|I| = 2 \sum i_k + l$) such that \medskip

$\bullet \;$ the resulting composites 

\[
\xymatrix{S^{-|I|-1} \ar[d] \ar[r]^{S^{-1}\zeta_I} & S^{-1}L_{K(n)} \ar[r] & [L_{n-1},L_n] \\ 
S^{-|I|-1}L_{n-1-i_l} \ar@{.>}[urr]^{S^{-1}\x_I} }
\]

\noindent
factor as shown, yielding an equivalence

\[
\xymatrix{
\vee S^{-1} \x_I : \vee_I \; S^{-|I|-1} L_{n-1-i_l} \ar[r]^\cong & [L_{n-1},L_n] \;. }
\]
\bigskip

\section{Generating functions}\bigskip

\noindent
The vector space dual

\[
X \mapsto \sL^{*}_n(X) := [L_0,L_n \wedge X]^D_{-*} := \Hom([L_0,L_n \wedge X]_{-*},\Q)
\]

\noindent
of $[L_0,L_n \wedge X]$ defines a cohomology theory which, following Hovey, records the 
divisible groups in the homotopy of $L_n$. Let

\[
\sL_n(T) := \sum_s [L_0,L_n]^D_{-s} \cdot T^s \;,
\]

\noindent
where $\Q_p$-vector spaces have been identified with their dimensions, eg

\[
\sL_0(T) = 1, \; \sL_1(T) = T^2, \; \sL_2(T) = 2T^4 + T^5 \;.
\] 
\bigskip

\noindent
{\bf Proposition:} The splitting conjecture implies that

\[
\sum_{n \geq 0} \sL_n(T) \; u^n = (1 - \sum_{k \geq 1} \epsilon_{k-1}(T) \; (uT^2)^k)^{-1}
\]

\noindent
where 

\[
\epsilon_k(T) = \prod_{0 \leq i \leq k-1}(1 + T^{2i+1})
\]

\noindent
is the Poincar\'e series for $E^*(k) := H_c^*(S_k,\Q_p)$ (with $\epsilon_0 = 1$). \bigskip

\noindent
{\bf Proof:} It follows from the conjecture that

\[
[L_0,L_n] = [L_0,[L_{n-1},L_n]]
\]

\noindent
[2 Prop 5.1], so when $n \geq 1$,

\[
\sL_n(T)  = \sum [L_0,L_{n-1-i_l}]^D_{|I|+1-s} \cdot T^s = 
\sum \sL_{n-1-i_l}(T) \cdot T^{|I|+1} \;. 
\]

\noindent
This can be rearranged as 

\[
\sum_{0 \leq k \leq n-1} \sL_{n-1-k}(T) \cdot \sum_{i_l=k}  T^{|I|+1} \;;
\]

\noindent
but the right-most term can be written as the sum

\[
T \cdot T^{2k+1} \sum T^{|\tilde{I}|} = T^{2(k+1)} \epsilon_k(T)
\]

\noindent
over sequences $\tilde{I}$ of the form $\{0 \leq i_1 < \dots < k\}$.
\bigskip

\noindent
We thus have the recursion relation

\[
\tL_n(T) := T^{-2n}L_n(T) = \sum_{0 \leq k \leq n-1} \epsilon_{n-k-1}(T) \cdot \tL_k(T)
\]

\noindent
and if $\epsilon(u) := 1 + \sum_{i \geq 0} \epsilon_i u^{i+1}$ then

\[
\tL := \sum_{n \geq 0} \tL_n(T) \; u^n \; = \sum_{n-1 \geq k \geq 0,n \geq 0} 
\epsilon_{n-k-1}(T) \cdot \tL_k(T)
\]
\[
= \sum_{n-1 \geq k \geq 0, n \geq 0} \epsilon_{n-k-1} u^{n-k} \cdot \tL_k \:u^k \; = \;
1 + (\epsilon -1) \: \tL \;,
\]

\noindent
so $\tL = (2 - \epsilon)^{-1}$. Replacing $u \mapsto uT^2$, we have

\[
\sum_{n \geq 0} \sL_n(T) \; u^n = (1 - \sum_{k \geq 1} \epsilon_{k-1}(T) \; (uT^2)^k)^{-1}
\]

\noindent
as asserted. $\Box$ 
\bigskip

\noindent
For example, if we specialize $T \to 1$ then $\epsilon_k \to 2^k$ and

\[
\epsilon \to 1 + \sum_{0 \leq k}2^ku^{k+1} = \frac{1-u}{1-2u} \;,
\]

\noindent
so 

\[
\tL \to \frac{1-2u}{1-3u} = 1 + \sum_{1 \leq k}3^{k-1}u^k \;,
\]

\noindent
ie $\sL_n(1) \to 3^{n-1}$ if $n \geq 1$ [H \S 5]. Similarly, we recover

\[
\sL_3(1) = T^6(1 + 2\epsilon_1 + \epsilon_2) = 4T^6 + 3T^7 + T^9 + T^{10}\;.
\] \bigskip

\noindent
{\bf Corollary:} The splitting conjecture implies that
\[
\sL_n(T) = \sum {n \choose \Sigma n_i} \; \epsilon_{i-1}(T)^{n_i} \cdot T^{2r} \;,
\]
(indexed by partitions  $n = n_1 + \cdots + n_r$ of $n$ with $r$ parts, with
\[
{n \choose \Sigma n_i} = \frac{n!}{\prod n_i!}
\]
denoting the multinomial coefficient). Equivalently, the cohomology theory $\sL_n$ is 
represented by the $p$-adic rationalization of the spectrum
\[
\bigvee S^{2r} (\prod \U(i-1)^{n_i})_+
\]
(indexed now by {\bf ordered} partitions of $n$), with $\U(k)$ the group of $k \times k$ 
unitary matrices. \bigskip

\noindent
{\bf Proof:} By the multinomial theorem. $\Box$ \bigskip

\noindent
{\bf Remarks:} \bigskip

\noindent
1) This corollary is a weakening of the original conjecture; but just because it is a 
question about cohomology theories in  characteristic zero doesn't seem to make it 
more accessible.\bigskip

\noindent
2) Mark Behrens asks if or how Gross-Hopkins duality fits into this story. The decomposition
above suggests the existence of a (noncommutative) product
\[
\sL_n \wedge \sL_m \to \sL_{n+m}
\]
and the proposition looks like it might have a reformulation in terms of some kind of Koszul 
duality; but I can't imagine what that might be \dots 
\bigskip

\bibliographystyle{amsplain}

\end{document}